\documentclass{amsart}

\usepackage{graphicx}
\usepackage{pst-all}      
\usepackage{url}

\newtheorem{theorem}{Theorem}

\newtheorem{conjecture}{Conjecture}

\begin{document}

\title{Irreducible triangulations of low genus surfaces}


\author{Thom Sulanke}
\address{Department of Physics, Indiana University, Bloomington, 
	Indiana 47405}
\email{tsulanke@indiana.edu}

\date{\today}

\begin{abstract}
The complete sets of irreducible triangulations are known for 
the orientable surfaces with genus of 0, 1, or 2 and 
for the nonorientable surfaces with genus of 1, 2, 3, or 4.
By examining these sets we determine some of the properties of these 
irreducible triangulations.
\end{abstract}

\maketitle

\section{Introduction}
\label{intro}

The irreducible triangulations of a surface provide a basis for obtaining all
the triangulations of that surface.
We can sequentially contract edges of a triangulation until an 
irreducible triangulation is produced.
Reversing this sequence we can produce any triangulation of a surface 
with a sequence of vertex splittings
starting with an irreducible triangulation.
Thus all the triangulations of a surface can be generated from the irreducible 
triangulations of that surface by vertex splittings.
The irreducible triangulations of a surface can be used to actually generate 
\cite{surftri} the triangulations of the surface.

Irreducible triangulations can also be used to check properties which are 
preserved by vertex splitting.
For example, let $\mathcal{P}$ be a property possessed by some of the 
triangulations of the surface $S$ which is 
preserved by vertex splitting such as ``contains a cycle which separates 
the surface $S$''.
If every irreducible triangulations of $S$ possesses $\mathcal{P}$ 
then every triangulations of $S$ possesses $\mathcal{P}$.
Conversely, if there is a counterexample to ``all triangulations of $S$ 
possess $\mathcal{P}$'' then there is a counterexample among the irreducible 
triangulation of $S$.

For any fixed surface the number of irreducible triangulations is finite
\cite{MR1021367}.
Irreducible triangulations have been determined and displayed by a number of 
authors:
the single irreducible triangulation of the sphere ($S_0$) by 
Steinitz and Rademacher \cite{StRa};
the two irreducible triangulations of the projective plane or the cross 
surface ($N_1$) by Barnette \cite{MR84f:57009}; 
the $21$ irreducible triangulations of the torus ($S_1$) by 
Lawrencenko \cite{MR914777});
and the $29$ irreducible triangulations of the Klein bottle ($N_2$) 
by Lawrencenko and Negami \cite{MR98h:05067} and 
Sulanke \cite{math.CO/0407008}.
The irreducible triangulations of the double torus ($S_2$), 
the triple cross surface ($N_3$), and the quadruple cross surface ($N_4$)
have been generated by the author using a computer program \cite{gentriang}.

\section{Definitions}
\label{definitions}

A {\em triangulation\/} of a closed surface 
is a simple graph embedded in the surface
such that each face is a triangle and any two faces share at most one edge.   

In a triangulation $T$ let
$abc$ and $acd$ be two faces which have $ac$ as a common edge.
The {\em contraction\/} of $ac$ is obtained by deleting $ac$,
identifying vertices $a$ and $c$, 
removing one of the multiple edges $ab$ or $cb$,
and removing one of the multiple edges $ad$ or $cd$.
The edge $ac$ of a triangulation $T$ is {\em contractible\/} if the 
contraction of $ac$ yields another triangulation of the surface in
which $T$ is embedded. 
If the edge $ac$ is contained in a 3-cycle other than the two 
which bound the faces which share it then its contraction 
would produce multiple edges.
Thus, for a triangulation $T$, not $K_4$ embedded in the sphere, 
an edge of $T$ is not contractible
if and only if that edge is contained in at least three 3-cycles. 
A triangulation is said to be {\em irreducible\/} 
if it has no contractible edges. 

The operation of {\em splitting} a vertex is the reverse of contracting an 
edge.
In a triangulation let $ab$ and $ac$ be two distinct edges.
The {\em splitting} of the vertex $a$ (along the edges $ab$ and $ac$) 
is obtained by
creating a new vertex $a'$, three new edges $a'a$, $a'b$, and $a'c$, and two
new faces $a'ab$ and $a'ac$.
The triangulation obtained by splitting a vertex is embedded in the same 
surface as the original triangulation.

We denote the orientable surface with genus $g$, the sphere with $g$ handles
attached, as $S_g$ and the 
nonorientable surface with genus $g$, the sphere with $g$ crosscaps attached,
as $N_g$.
Define the {\em Euler genus} of the surface $S$ to be $2 - \chi(S)$.
For orientable surfaces the Euler genus is twice the genus and for
nonorientable surfaces the Euler genus is the same as the genus.

\section{Generating irreducible triangulations}

The author has recently developed an algorithm \cite{gentriang} 
for generating irreducible 
triangulations of a surface by using the irreducible triangulations of
other surfaces with smaller Euler genera.
This algorithm was implemented as a computer program.
The irreducible triangulations of $S_2$, $N_3$, and $N_4$ were generated
and are available as computer files \cite{surftri}.

Before we briefly describe the algorithm used to generate irreducible 
triangulations
we examine how an irreducible triangulation can be reduced to an irreducible
triangulation with a lower genus.
For simplicity we only consider orientable surfaces here.
Let $T$ be an irreducible triangulation of $S_g$ with $g>0$.
Every edge of $T$ is on a 3-cycle which is not a face.
Many of these 3-cycles do not separate $S_g$ into two components.
Pick one of these nonseparating 3-cycles.
Cut $T$ along this 3-cycle thereby cutting one of the handles of $S_g$.
Cap the resulting two holes with new
triangular faces to produce a new triangulation $T'$ of $S_{g-1}$. 
Contract contractible edges until an irreducible 
triangulation of $S_{g-1}$ is obtained.

To generate an irreducible triangulation of $S_g$ we reverse these steps
in effect ``growing a handle''.
Start with an irreducible triangulation of $S_{g-1}$.  
Split vertices checking each new triangulation to see if it
can be used to form an irreducible triangulation of $S_g$.
The final step is the reverse of the cut and cap described above.
Remove two faces and join the resulting boundary cycles
in such a way that the resulting triangulation is still orientable.

An irreducible triangulation of $N_g$ can be generated in a similar way by
``growing a handle or a crosshandle''.
Start with an irreducible triangulation of $S_{g/2-1}$ or $N_{g-2}$
and split vertices.
In the final step we remove two faces and join the resulting boundary cycles
in such a way that the resulting triangulation is nonorientable.

We can also ``grow a crosscap'' to generate an irreducible triangulation of
$N_g$ starting with an irreducible
triangulation of $S_{(g-1)/2}$ or $N_{g-1}$.
As new triangulations are produced by edge splitting we check for vertices with
degree 6.
When we remove a vertex with degree 6 and its incident faces 
a hole with a 6-cycle as a boundary is produced.
By identifying the 3 pairs of opposite vertices on this 6-cycle we check if
the result is an irreducible triangulation.

\section{Counts}

Due to the large number of irreducible triangulations of 
$S_2$, $N_3$, and $N_4$
the irreducible triangulations are not be displayed here but some of their 
properties are presented.
For comparison we also include similar properties for $S_0$, $S_1$, $N_1$, and
$N_2$.

Table \ref{irrtri} shows for each surface the number of irreducible 
triangulations having a given number of vertices.

\begin{table}
\centering
\begin{tabular}{r|r r r r r r}
Vertices & $S_1$ & $S_2$ & $N_1$ & $N_2$ & $N_3$ & $N_4$ \\
\hline
6     &    &        & 1 &    &      &         \\
7     &  1 &        & 1 &    &      &         \\
8     &  4 &        &   &  6 &      &         \\
9     & 15 &        &   & 19 &  133 &      37 \\
10    &  1 &    865 &   &  2 & 2521 &   10347 \\
11    &    &  26276 &   &  2 & 4638 &  370170 \\
12    &    & 117047 &   &    & 1320 & 1891557 \\
13    &    & 159205 &   &    &  946 & 2067817 \\
14    &    &  54527 &   &    &   93 &  956967 \\
15    &    &  38195 &   &    &   50 &  700733 \\
16    &    &    664 &   &    &    7 &  186999 \\
17    &    &      5 &   &    &      &   89036 \\
18    &    &        &   &    &      &   19427 \\
19    &    &        &   &    &      &    3975 \\
20    &    &        &   &    &      &     832 \\
21    &    &        &   &    &      &      79 \\
22    &    &        &   &    &      &       6 \\
\hline
Total & 21 & 396784 & 2 & 29 & 9708 & 6297982 \\
\end{tabular}
\vspace{.1in}
\caption{Irreducible triangulation by vertices}
\label{irrtri}
\end{table}

\section{Noncontractible separating cycles}
\label{nsc}

Let $v_1 v_2 \ldots v_n$ be
an n-cycle in a graph embedded on the surface $S$ and let $C$ be the closed
curve which is the embedding of $v_1 v_2 \ldots v_n$ in $S$.
$v_1 v_2 \ldots v_n$ is {\em separating\/} if $S - C$ is disconnected.
$v_1 v_2 \ldots v_n$ is {\em contractible\/} if a component of $S - C$ is 
a 2-cell, otherwise, it is {\em noncontractible\/}.
This definition of a contractible cycle should not be confused with the 
definition of a contractible edge given earlier.
Necessary conditions for the existence of a
{\em noncontractible separating cycle} or {\em NSC} have been studied
\cite{MR1996201} \cite{MR1118042} \cite{MR1224715}.
An NSC separates a surface into two components neither of which is a 2-cell.
Thus a surface having an NSC must have genus greater than 1.

The existence of an NSC in a triangulation and
the genera of the separated surfaces  
are preserved by vertex splitting.
Thus if every irreducible triangulation of a surface has an NSC then every
triangulation of that surface has an NSC.

Barnette conjectured that every triangulation of a
surface with genus greater than 1 has an NSC.
Lawrencenko and Negami \cite{MR98h:05067} showed that every irreducible 
triangulation of $N_2$ has an NSC and thus every triangulation of $N_2$ 
has an NSC.
Ellingham, Zha, and Jennings \cite{dana} have shown (without any reference to 
irreducible triangulations) that every triangulation of $S_2$ has an NSC.

By checking that each irreducible triangulation of $S_2$, $N_2$, $N_3$, and 
$N_4$ has an NSC we have the following result which is new only for $N_3$, and 
$N_4$.

\begin{theorem}
\label{nscs2n2n3n4}
Every triangulation of $S_2$, $N_2$, $N_3$, or $N_4$ has an NSC.
\end{theorem}

Similarly, if an NSC separates a surface with Euler genus $g$ into two surfaces
with Euler genera $h$ 
and $g-h$ then any triangulation obtained by vertex splitting of this 
triangulation has an NSC which separates the surface into two surfaces with 
Euler genera $h$ and $g-h$.

Thomassen conjectured (\cite{MR1844449} page 167) that given a triangulation 
of an orientable surface with genus $g$ and an integer $h$ such that 
$1 \leq h < g$, 
then the triangulation must contain an NSC such that the two surfaces 
separated by the NSC (after capping the holes with disks) 
have genera $h$ and $g-h$, respectively.
This conjecture is equivalent to Barnette's conjecture for $S_2$ (and $S_3$)
but we can make a similar conjecture for nonorientable surfaces.

\begin{conjecture}
\label{nonornscallgenus}
Given a triangulation of a nonorientable surface with Euler genus $g$ and 
an integer $h$ such that $1 \leq h < g$, 
then the triangulation must contain an NSC such that the two surfaces 
separated by the NSC have Euler genera $h$ and $g-h$, respectively.
\end{conjecture}

By checking the irreducible triangulations of $N_4$ we have the following.

\begin{theorem}
\label{nonornscallgenusn4}
Every triangulation of $N_4$ has an NSC which separates the surface into two
surfaces each with Euler genus 2.
Every triangulation of $N_4$ has an NSC which separates the surface into two
surfaces with Euler genera 1 and 3, respectively.
\end{theorem}

From Theorems~\ref{nscs2n2n3n4} and~\ref{nonornscallgenusn4} it follows that
Conjecture~\ref{nonornscallgenus} is true for $N_2$, $N_3$, and $N_4$.
Conjecture~\ref{nonornscallgenus} and Theorem~\ref{nonornscallgenusn4} do not 
specify the orientability of the separated surfaces.
For example,
there are triangulations of $N_3$ which do not have an NSC which separates the
surface into $N_1$ and $N_2$.
Such an example can be constructed using any irreducible triangulation of 
$N_1$ and any irreducible triangulation of $S_1$.  
Remove a face from each of these two irreducible triangulations and identify 
the resulting boundaries.
Let $C_1$ be the closed curve in $N_3$ where the two surfaces were joined.
Assume there is an NSC which separates $N_3$ into $N_1$ and $N_2$.
Let $C_2$ be the closed curve in $N_3$ corresponding to this NSC.
Due to the topology of $N_3$ the curves $C_1$ and $C_2$ must cross at least 
four times.
But this contradicts the fact that the cycle corresponding to $C_1$ has 
length 3.
Similarly, there are also triangulations of $N_3$ which do not have an NSC 
which separates the surface into $N_1$ and $S_1$.

For $N_3$, there are 9184 irreducible triangulations which have an NSC which 
separates the surface into $N_1$ and $N_2$ and 
there are 8533 irreducible triangulations which have an NSC which 
separates the surface into $N_1$ and $S_1$.

For $N_4$, there are 6062415 irreducible triangulations which have an NSC 
which separates the surface into $N_2$ and $N_2$ and 
there are 5971981 irreducible triangulations which have an NSC which 
separates the surface into $N_2$ and $S_1$.

The {\em edge-width\/} of a triangulation is the length of the shortest 
NSC in the triangulation.
Tables \ref{minnscs2} through \ref{minnscn4} show the number of irreducible 
triangulations for a given number of vertices and a given value of the 
edge-width.

\begin{table}
\centering
\begin{tabular}{r|r r r r r r}
          & \multicolumn{6}{c}{Edge-width} \\
Vertices  & 3 &     4 &     5 &     6 &     7 &   8 \\
\hline
10    &       &     2 &    51 &   681 &   130 &   1 \\
11    &     2 &    58 &  2249 & 16138 &  7818 &  11 \\
12    &    25 &  1516 & 20507 & 72001 & 22877 & 121 \\
13    &   710 & 13004 & 50814 & 78059 & 16609 &   9 \\
14    &  8130 & 30555 & 12308 &  3328 &   205 &   1 \\
15    & 36794 &  1395 &     3 &     1 &     2 &     \\
16    &   661 &     3 &       &       &       &     \\
17    &     5 &       &       &       &       &     \\
\hline
\end{tabular}
\vspace{.1in}
\caption{Irreducible triangulation of $S_2$ by vertices and edge-width}
\label{minnscs2}
\end{table}


\begin{table}
\centering
\begin{tabular}{r|r r r r r r}
          & \multicolumn{6}{c}{Edge-width} \\
Vertices  & 3 & 4 & 5 & 6 & & \\
\hline
8     &   & 1 & 5 &    &  &   \\
9     & 1 & 5 & 2 & 11 &  &   \\
10    & 1 & 1 &   &    &  &   \\
11    & 2 &   &   &    &  &   \\
\hline
\end{tabular}
\vspace{.1in}
\caption{Irreducible triangulation of $N_2$ by vertices and edge-width}
\label{minnscn2}
\end{table}

\begin{table}
\centering
\begin{tabular}{r|r r r r r r}
          & \multicolumn{6}{c}{Edge-width} \\
Vertices  & 3 & 4 & 5 & 6 & & \\
\hline
9     &     &    1 &  119 &   13 &  &   \\
10    &   1 &  140 & 1862 &  518 &  &   \\
11    &  72 & 1248 & 1558 & 1760 &  &   \\
12    & 502 &  811 &    4 &    3 &  &   \\
13    & 912 &   34 &      &      &  &   \\
14    &  93 &      &      &      &  &   \\
15    &  50 &      &      &      &  &   \\
16    &   7 &      &      &      &  &   \\
\hline
\end{tabular}
\vspace{.1in}
\caption{Irreducible triangulation of $N_3$ by vertices and edge-width}
\label{minnscn3}
\end{table}

\begin{table}
\centering
\begin{tabular}{r|r r r r r r}
          & \multicolumn{6}{c}{Edge-width} \\
Vertices  & 3 & 4 & 5 & 6 & 7 & 8 \\
\hline
9     &        &        &     17 &     20 &  \\
10    &        &        &   5028 &   5222 &    97  \\
11    &        &   4503 & 209623 & 150994 &  5050  \\
12    &   2499 & 161502 & 983249 & 717138 & 27169  \\
13    &  76309 & 704856 & 698076 & 566851 & 21723 & 2  \\
14    & 396148 & 519038 &  36649 &   5066 &    66  \\
15    & 633195 &  67538  \\
16    & 181884 &   5115  \\
17    &  88799 &    237  \\
18    &  19427  \\
19    &   3975  \\
20    &   832   \\
21    &    79   \\
22    &     6   \\
\hline
\end{tabular}
\vspace{.1in}
\caption{Irreducible triangulation of $N_4$ by vertices and edge-width}
\label{minnscn4}
\end{table}


\section{Nonseparating cycles}

Every cycle of a triangulation of $S_0$ separates and we exclude such 
triangulations in this section.
In \cite{gentriang} it is shown that for every vertex of an irreducible 
triangulation 
there are at least two nonseparating 3-cycles containing that vertex.
Thus every irreducible triangulation has a nonseparating cycle and, therefore,
every triangulation has a nonseparating cycle.
For triangulations of orientable surfaces the only topological type
of a nonseparating cycle is one which cuts a handle.

Let $S$ be a triangulated nonorientable surface with a nonseparating cycle. 
Let $C$ be the closed curve which the embedding of that cycle in $S$.
The cycle is {\em one-sided} if the neighborhood of $C$ in $S$ is 
homeomorphic to a M\"{o}bius band, otherwise the cycle is {\em two-sided}.
The cycle is {\em orientable-leaving} if $S-C$ is orientable, otherwise
the cycle is {\em nonorientable-leaving}.
For triangulations of nonorientable surfaces there are four possible 
topological types of nonseparating cycles depending on whether it is one- or 
two-sided and whether it is orientable- or nonorientable-leaving.
At most three of the these types of nonseparating cycles can occur for a fixed
nonorientable surface since orientable surfaces have an even Euler genus.

By checking the irreducible triangulations of $N_1$, $N_2$, $N_3$, and $N_4$
we have the following theorem.

\begin{theorem}
\label{nonsep}
Every triangulation of $N_1$ has
a nonseparating cycle which is one-sided and orientable-leaving.

Every triangulation of $N_2$ has 
a nonseparating cycle which is one-sided and nonorientable-leaving and 
a nonseparating cycle which is two-sided and orientable-leaving.

Every triangulation of $N_3$ has 
a nonseparating cycle which is one-sided and orientable-leaving; 
a nonseparating cycle which is one-sided and nonorientable-leaving; and 
a nonseparating cycle which is two-sided and nonorientable-leaving.

Every triangulation of $N_4$ has 
a nonseparating cycle which is one-sided and nonorientable-leaving; 
a nonseparating cycle which is two-sided and orientable-leaving; and
a nonseparating cycle which is two-sided and nonorientable-leaving.
\end{theorem}

\begin{conjecture}
If $g>=3$ is odd then every triangulation of $N_g$ has 
a nonseparating cycle which is one-sided and orientable-leaving; 
a nonseparating cycle which is one-sided and nonorientable-leaving; and 
a nonseparating cycle which is two-sided and nonorientable-leaving.

If $g>=4$ is even then every triangulation of $N_g$ has 
a nonseparating cycle which is one-sided and nonorientable-leaving; 
a nonseparating cycle which is two-sided and orientable-leaving; and
a nonseparating cycle which is two-sided and nonorientable-leaving.
\end{conjecture}

\section{Maximal irreducible triangulations}

Define $V_{max}(S)$ to be the maximum number of vertices in an irreducible
triangulation of $S$.
From Tables \ref{minnscs2} through \ref{minnscn4} we see that if an irreducible
triangulation $T$ of one of these surfaces $S$ has $|V(T)| = V_{max}(S)$ 
then $T$ has an NSC of length 3.
For $S_2$ this confirms a conjecture of Negami \cite{negamiconj}.
This suggests that maximal irreducible triangulations are made up of other
triangulations joined at a single face of each.

We can use the following construction to obtain large irreducible 
triangulations which give a lower bound for $V_{max}(S)$.
This construction is similar to the one given by Nakamoto and Ota 
\cite{MR1348564} and the lower bound which it provides is a slight improvement.

For $N_1$ there is only one maximal irreducible triangulation $M$ which has 7 
vertices.
If we take two copies of $M$ and remove one face from each we can join them
at the boundaries of these faces to obtain a triangulation of $N_2$. 
This triangulation has 11 vertices and is irreducible.

For each $g>2$ we will construct a base triangulation $B_g$ of $S_0$ which when
joined with $g$ copies of $M$ will produce an irreducible triangulation of 
$N_g$.

\begin{figure}
\centering
\psset{unit=.015\textwidth}
\parbox{.49\textwidth}{%
\centering
\ttfamily
\begin{pspicture}(-5,-5)(25,20)
\pspolygon*[linecolor=lightgray](0,0)(10,17)(10,6)
\pspolygon*[linecolor=lightgray](20,0)(10,17)(10,6)
\psline[linewidth=2pt]{C-C}(0,0)(20,0)
\psline[linewidth=2pt]{C-C}(0,0)(10,17)
\psline[linewidth=2pt]{C-C}(0,0)(10,6)
\psline[linewidth=2pt]{C-C}(20,0)(10,17)
\psline[linewidth=2pt]{C-C}(20,0)(10,6)
\psline[linewidth=2pt]{C-C}(10,17)(10,6)
\rput*(10,-3){$B_3$}
\end{pspicture}}
\parbox{.49\textwidth}{%
\centering
\ttfamily
\begin{pspicture}(-5,-5)(25,20)
\pspolygon*[linecolor=lightgray](10,17)(7,7)(13,7)
\pspolygon*[linecolor=lightgray](0,0)(10,3)(7,7)
\pspolygon*[linecolor=lightgray](20,0)(10,3)(13,7)
\psline[linewidth=2pt]{C-C}(0,0)(20,0)
\psline[linewidth=2pt]{C-C}(0,0)(10,17)
\psline[linewidth=2pt]{C-C}(0,0)(10,3)
\psline[linewidth=2pt]{C-C}(0,0)(7,7)
\psline[linewidth=2pt]{C-C}(20,0)(10,17)
\psline[linewidth=2pt]{C-C}(20,0)(10,3)
\psline[linewidth=2pt]{C-C}(20,0)(13,7)
\psline[linewidth=2pt]{C-C}(10,17)(7,7)
\psline[linewidth=2pt]{C-C}(10,17)(13,7)
\psline[linewidth=2pt]{C-C}(10,3)(7,7)
\psline[linewidth=2pt]{C-C}(10,3)(13,7)
\psline[linewidth=2pt]{C-C}(7,7)(13,7)
\rput*(10,-3){$B_4$}
\end{pspicture}}

\caption{Base triangulations for constructing large irreducible triangulations}
\label{bases}
\end{figure}
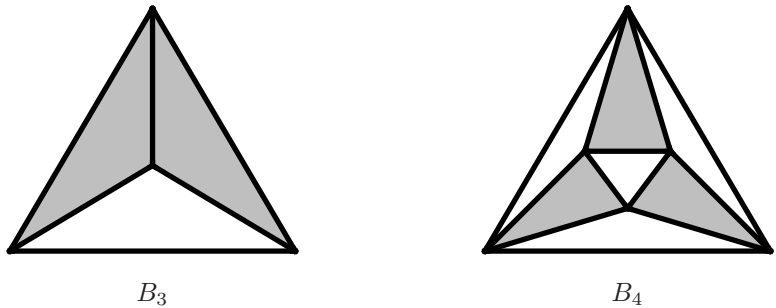

The left side of Figure~\ref{bases} shows $B_3$ which is
a triangulation of $S_0$ from which three faces
(the shaded faces and the outside face) have been removed.
Every edge of $B_3$ is on a removed face.
If we join three punctured copies of $M$ at these faces we get a
triangulation of $N_3$.
This triangulation is irreducible since each edge of $B_3$ is now on at least
three 3-cycles.

The right side of Figure~\ref{bases} shows $B_4$ 
which is a triangulation of $S_0$ from which four faces have been removed.  
Again every edge of $B_4$ is on a removed face.
When we join four punctured copies of $M$ at the removed faces we obtain an
irreducible triangulation of $N_4$.

If we take two copies of $B_4$ and join them at removed faces then 
we obtain $B_6$ which is a triangulation of $S_0$ with 6 faces removed.
Every edge of $B_6$ is either on a removed face or on at least three 3-cycles.
Joining six punctured copies of $M$ we obtain an
irreducible triangulation of $N_6$.
We can repeat this construction to obtain $B_g$ for even $g>2$.
$|V(B_4)| = 6$ and each additional copy of $B_4$ adds 3 vertices such that
$|V(B_g)| = 3g/2$.
Each copy of $M$ adds 4 vertices.
Thus for even $g>2$ the number of vertices in the constructed irreducible 
triangulation of $N_g$ is $11g/2$.

To obtain $B_g$ for odd $g$ we join $B_3$ to $B_{g-1}$.
Then $|V(B_g)| = 3(g-1)/2 + 1 = 3g/2 - 1/2$
and the number of vertices in the constructed irreducible 
triangulation of $N_g$ is $11g/2-1/2$.

Thus for any $g$ we have 
\[ V_{max}(N_g) \geq \lfloor{11g/2}\rfloor \]

For $S_1$ the only maximal irreducible triangulation has 10 vertices. 
Repeating the above construction with this triangulation as $M$ we obtain 
\[ V_{max}(S_g) \geq \lfloor{17g/2}\rfloor \]

In the above construction the triangulation $M$ does not need to be 
irreducible.
Any edge of the removed face may be contractible and the resulting 
triangulation would still be irreducible.

A triangulation is {\em almost irreducible\/} if it is not irreducible and it 
has a face which is incident on all the contractible edges.
If $M$ is almost irreducible then the construction still produces an 
irreducible triangulation.
However, there are no almost irreducible triangulations of $N_1$ 
\cite{MR98h:05067} and
there are no almost irreducible triangulations $T$ of $S_1$ for which 
$|V(T)| > V_{max}(S_1)$.
There are 8 almost irreducible triangulations of $S_1$ but none have more than
9 vertices \cite{almost}.

\section{Pseudo-minimal triangulations}

Two triangulations $T$ and $T'$ of
a surface are {\em equivalent\/} 
if there is a isomorphism $h$
with $h(T)=T'$.
That is, if $a$, $b$, and $c$ are vertices of $T$ then
$ab$ is an edge of $T$ if and only if $h(a)h(b)$ is an edge
of $T'$ and a face of $T$ is bounded by the cycle $abc$
if and only if a face of $T'$ is bounded by the cycle
$h(a)h(b)h(c)$.

Let $ac$ be an edge in a triangulation $T$ and
$abc$ and $acd$ be the two faces which have $ac$ as a common edge.
The {\em diagonal flip\/} of $ac$ is obtained by deleting $ac$,
adding edge $bd$, deleting the faces $abc$ and $acd$, and adding the faces
$abd$ and $bcd$.
An edge $ac$ of a triangulation $T$ is {\em flippable\/} if the 
diagonal flip of $ac$ yields another triangulation of the surface in
which $T$ is embedded.  
Thus the edge $ac$ is flippable if there is not already an edge $bd$.
Two triangulations are {\em equivalent under diagonal flips\/}
if one is equivalent to a triangulation obtained from the other 
by a sequence of diagonal flips.

The number of vertices of an irreducible triangulation 
can not be reduced by edge contraction.
Negami \cite{MR95m:05091} defines a type of triangulation for which 
the number of vertices can not be reduced by a combination of diagonal flips 
and edge contractions. 
An irreducible triangulation is said to be {\em pseudo-minimal\/} 
if it is equivalent under diagonal flips only to irreducible triangulations.

A triangulation is said to be {\em minimal\/} if there are no triangulations 
of the same surface with fewer vertices.  
It is clear that such a triangulation is also pseudo-minimal.  
The number of vertices in a minimal triangulation for
nonorientable surfaces was determined by Ringel \cite{MR17:774b} and for 
orientable surfaces by
Jungerman and Ringel \cite{MR82b:57012}.  
It is given for all surface except $N_2$, $N_3$, and $S_2$ by the formula:
\[ V_{min}(S) = \left\lceil{\frac{7+\sqrt{49-24\chi(S)}}{2}}\right\rceil \]
For the three exceptions the value is one more than the value given by the 
formula:
$V_{min}(N_2)=8$, $V_{min}(N_3)=9$, and $V_{min}(S_2)=10$.

Let $N(S)$ be the minimum value such that two triangulations $T$ and $T'$
of $S$ are equivalent under diagonal flips 
if $|V(T)| = |V(T')| \geq N(S)$.
Negami \cite{MR95m:05091} has
shown that such a finite value exists for any $S$.

$N(S_0) = V_{min}(S_0) = 4$, $N(S_1) = V_{min}(S_1) = 7$, 
$N(N_1) = V_{min}(N_1) = 6$, and $N(N_2) = V_{min}(N_2) = 8$ 
are known \cite{wagner} 
\cite{MR46:8878} \cite{MR91g:05038}.

Checking the irreducible triangulations generated for $S_2$
we have determined that the 865 minimal triangulations are the only 
pseudo-minimal triangulations and that these pseudo-minimal triangulations 
form one equivalence class under diagonal flips.
Thus $N(S_2) = V_{min}(S_2) = 10$ (\cite{MR95m:05091} \cite{examples}).
Similarly, the 133 minimal triangulations $N_3$ are the only 
pseudo-minimal triangulations and they 
form one equivalence class under diagonal flips.
Thus $N(N_3) = V_{min}(N_3) = 9$.

The situation for $N_4$ is more interesting.  
The 37 minimal triangulations are the only pseudo-minimal triangulations.
However, these pseudo-minimal triangulations are partitioned into 
three equivalence 
classes under diagonal flips \cite{examples} with cardinality 32, 3, and 2.
Using this complete list of pseudo-minimal triangulations of $N_4$ it is
possible to show \cite{examples} that $N(N_4) = V_{min}(N_4) + 1 = 10$.

Suppose for a surface $S$ there exist at least two inequivalent minimal 
triangulations which have no flippable edges.
Then $N(S) > V_{min}(S)$.
There are an infinite number of surfaces
which have many inequivalent triangular embeddings of complete graphs
\cite{MR2000i:05055} \cite{MR2064871} \cite{MR2001k:05058}.
A triangular embeddings of complete graph is minimal and a complete graph
has no flippable edges.
Therefore there are an infinite number of surfaces $S$ for which
$N(S) > V_{min}(S)$.
For each of the surfaces $S_g, 3<=g<=15$ and $N_g, 5<=g<=30$, the author has
found, using random computer searching \cite{examples},
a pair of minimal triangulations which are inequivalent
under diagonal flips.
The existence of these pairs shows that if $3<=g<=15$ then
$ N(S_g) > V_{min}(S_g) $
and that if $4<=g<=30$ then
$ N(N_g) > V_{min}(N_g) $.

\begin{conjecture}
The only surfaces $S$ for which  $N(S) = V_{min}(S)$ are 
$S_0$, $S_1$, $S_2$, $N_1$, $N_2$, and $N_3$.
\end{conjecture}



\bibliographystyle{amsplain}
\bibliography{properties}

\end{document}